\documentclass{article}

\usepackage[cp1251]{inputenc}
\usepackage[russian]{babel}
\usepackage{amsfonts}
\usepackage{amsmath}
\usepackage{amsfonts,amssymb}
\usepackage{amssymb}

\begin{document}

\def\eps{\varepsilon}
\def\W{W_2^{-1}[0,1]}

\centerline{\large\bf Экстремальные свойства первого собственного значения}
\centerline{\large\bf задачи Штурма-Лиувилля с краевыми условиями третьего типа.}

\bigskip
\centerline{\bf Карулина Е.С.}

\bigskip
\bigskip
Рассматривается задача Штурма---Лиувилля
\begin{gather} \label{Karulina_eq1}
y''+ qy+\lambda y=0,\\\label{Karulina_eq2}
y'(0)-k_0^2y(0)= y'(1)+k_1^2y(1)=0,
\end{gather}
\vskip -1ex\noindent
где $k_0,\,k_1 \in \mathbb{R}$, а функция $q$ принадлежит множеству
$$
A_\gamma = \{q\in L_1[0,1] : q(x)\geqslant 0,
\;\int_0^1 q^\gamma\,dx=1\}$$
при $\gamma\in \mathbb{R} \setminus \{0\}$.

\medskip
Пусть $m_{\gamma}=\inf\limits_{q\in A_{\gamma}}\lambda_{1}(q)$, $M_{\gamma}=\sup\limits_{q\in A_{\gamma}}\lambda_{1}(q)$, где $\lambda_{1}$ --- первое собственное значение задачи \eqref{Karulina_eq1}--\eqref{Karulina_eq2}.
Цель данной работы "--- найти значения \(m_{\gamma}\) и \(M_{\gamma}\) при некоторых значениях $\gamma$.

\bigskip
Значения \(m_{\gamma}\) и \(M_{\gamma}\) для уравнения
$y''+\lambda q(x) y = 0$
с условиями Дирихле были найдены Егоровым и Кондратьевым в 1996 г. (см., например, \cite{EK:1996}). Ими впервые была рассмотрена задача такого типа.

Значение
\(m_\gamma\) для задачи Дирихле с уравнением \eqref{Karulina_eq1}
было рассмотрено в 2003 г. Винокуровым и Садовничим при $\gamma\geq 1$ (см. \cite{VS:2003}).
В работах Ежак С.С. для этой задачи были получены значения \(m_{\gamma}\) и \(M_{\gamma}\) при всех значениях $\gamma$ (см., например, \cite{Ez:2015}).
Также в работе Владимирова А.А. в 2016 г. было уточнено значение \(m_\gamma\) при $\gamma\in[1/3,1/2)$ (см. \cite{Vl:2016}).

\bigskip
Задача $y''- qy + \lambda y=0$ с условиями \eqref{Karulina_eq2}, $q\in A_\gamma$ подробно рассмотрена в работах автора (см., например, \cite{Kar:Uniti}).

\medskip
В 2013 г. в \cite{K_V:Tatra} были найдены значения \(m_1\) и $M_1$ для задачи \eqref{Karulina_eq1}--\eqref{Karulina_eq2}.

\bigskip
В 2016 г. Ежак С.С. и Тельновой М.Ю. получены первые результаты для задачи Дирихле с уравнением \eqref{Karulina_eq1} и весовым интегральным условием (см., например, \cite{TE:2017}).

\bigskip
Основным результатом данной работы являются следующие теоремы:

\textbf{Теорема 1.}
{\it Если $\gamma<1$, то $m_\gamma=-\infty$.}

\textbf{Теорема 2.}
{\it Если $\gamma>1$, то $M_\gamma$ --- это первое собственное значение задачи $y''+\lambda y=0$, \eqref{Karulina_eq2}.}

\bigskip

Мы предполагаем, что все рассматриваемые пространства являются вещественными.

В данной статье мы расширяем класс допускаемых к рассмотрению потенциалов
с пространства \(L_1[0,1]\) до пространства \(W_2^{-1}[0,1]\) (см. \cite{Vl:2004}, \cite{K_V:Tatra}). Это пространство, в частности, содержит всевозможные
дельта-функции. Такое обобщение рассматриваемой задачи позволяет получить нужные нам оценки, а также доказать, что они достигаются
на потенциалах из расширенного класса.

Через \(W_2^{-1}[0,1]\) будем обозначать гильбертово пространство, являющееся
пополнением пространства \(L_2[0,1]\) по норме
\[
	\|y\|_{W_2^{-1}[0,1]}\rightleftharpoons\sup\limits_{\|z\|_{W_2^1[0,1]}=1}
		\int_0^1 yz\,dx.
\]
Если \(y\in W_2^{-1}[0,1]\), то через \(\int_0^1 yz\,dx\) мы обозначаем результат
применения линейного функционала \(y\) к функции \(z\in W_2^1[0,1]\):
\[
	\langle y,z\rangle\rightleftharpoons\lim\limits_{n\to\infty}
		\int_0^1 y_nz\,dx \quad (\text{где } y=\lim_{n\to\infty}y_n,
		\quad y_n\in L_2[0,1]).
\]

Пусть \(\Gamma_\gamma\) "--- замыкание в пространстве \(W_2^{-1}[0,1]\) множества \(A_\gamma\). Тогда
$$
m_{\gamma}=\inf\limits_{q\in \Gamma_{\gamma}}\lambda_{1}(q),\quad M_{\gamma}=\sup\limits_{q\in \Gamma_{\gamma}}\lambda_{1}(q).
$$

\medskip
Для $p\in \mathbb{R}$ и $y\in L_\infty[0,1]:y^{-1}\in L_\infty[0,1]$ определим величину
$$
\|y\|_p = \lim_{\scriptstyle{\vphantom{b}r\to p}\atop\scriptstyle{r\ne 0}}\left(\int_{0}^{1}|y|^r dx\right)^{1/r}.
$$

Из асимптотического соотношения
\begin{equation}\label{norm0}
	\|y\|_0 = \lim_{\scriptstyle{\vphantom{b}r\to 0}\atop\scriptstyle{r\ne 0}}\left(\int_0^1 |y|^r\,dx\right)^{1/r}=
\exp\left(\int_0^1 \ln{|y|}\,dx\right)
\end{equation}
следует, что при фиксированном $y$ величина $\|y\|_p$ зависит от $p$ непрерывно.

\bigskip
\textbf{Утверждение 0.}
{\it Если $p,r\in \mathbb{R}$, $p<r$, то $\|y\|_p\leq\|y\|_r$.}

\bigskip
\textbf{Доказательство утверждения 0.}
\medskip

Используем неравенство Гёльдера.

Пусть $p>0$, тогда
$$ {
\|y\|_p = \left(\int_{0}^{1}|y|^p dx\right)^{1/p} \leqslant
\left(\left(\int_{0}^{1}|y|^r dx\right)^{p/r} \right)^{1/p} = \|y\|_r
}.
$$

Пусть $r<0$, тогда
$$ {
\|y\|_r = \left(\int_{0}^{1}|y|^r dx\right)^{1/r} \geqslant
\left(\left(\int_{0}^{1}|y|^p dx\right)^{r/p} \right)^{1/r} = \|y\|_p
}.
$$

Из \eqref{norm0} и непрерывности величины $\|y\|_p$ следует, что данное неравенство выполняется и для произвольных $p,r\in \mathbb{R}$.


\bigskip
Обобщённая функция \(q\in W_2^{-1}[0,1]\) называется
\textit{неотрицательной}, если для любой неотрицательной функции \(y\in W_2^1[0,1]\) выполняется
неравенство \(\langle q,y\rangle\geq 0\).

\bigskip
\textbf{Утверждение 1.} 
{\it Пусть $\gamma\in(0,1)$, $\eps>0$, $\rho\in\W$, $\rho\ge 0$. Тогда найдется последовательность $\{q_n\}$ неотрицательных функций из $L_\infty[0,1]$ со свойствами $\|q_n\|_\gamma<\eps$ для достаточно больших $n$,
$\rho=\lim\limits_{n\to\infty}{q_n}$ в $\W$.}

\bigskip
\textbf{Доказательство утверждения 1.}
\medskip

Пусть $\zeta\in[0,1]$,
$$
q_n(x) = \left\{
\begin{array}{ll}
n, & x-(\zeta-1/n)^+\in(0,1/n), \\
0, & \mbox{ иначе}.
\end{array}
\right.
$$

Символом $a^+$ здесь и далее обозначается {\it положительная часть\/} числа $a$,
то есть число
$$
	a^+=\begin{cases}
a& \mbox{если } a>0,\cr
0&\mbox{иначе}.
\end{cases}
$$

Для любого $\eps>0$ при достаточно больших $n$ выполняется неравенство
$$
\|q_n\|_\gamma = n^\frac{\gamma-1}{\gamma}<\eps. 
$$

Докажем, что последовательность $\{q_n\}$ фундаментальна в $\W$ и что $\lim\limits_{n\to\infty}{q_n}=\delta_{\zeta}$\footnote{Здесь и далее символом \(\delta_{\zeta}\) мы обозначаем
дельта-функцию Дирака с носителем в точке \(\zeta\in [0,1]\).}. Поскольку любая функция из $\W$ может быть аппроксимирована линейной комбинацией дельта-функций, этого будет достаточно для доказательства утверждения 1. 

По определению нормы в $\W$, для любых $n,m\in \mathbb{N}$
$$
\|q_n-q_m\|_{W_2^{-1}[0,1]} = \sup\limits_{\|z\|_{W_2^1[0,1]}=1}
		\int_0^1 (q_n-q_m)z\,dx = \sup\limits_{\|z\|_{W_2^1[0,1]}=1}
		(z(\xi_z)-z(\eta_z)), 
$$
где $\left|\xi_z-\eta_z\right|<\max{\left(1/n,1/m\right)}$.

Используем неравенство Гельдера:
$$
|z(\xi_z)-z(\eta_z)| \leq \|z\|_{W_2^1[0,1]}\sqrt{|\xi_z-\eta_z|},
$$
откуда следует
$$
\|q_n-q_m\|_{W_2^{-1}[0,1]} \leq \sup\limits_{\|z\|_{W_2^1[0,1]}=1}
		\sqrt{|\xi_z-\eta_z|} < \varepsilon \mbox{ при } n,m>\frac{1}{\varepsilon}.
$$
Следовательно, последовательность $\{q_n\}$ фундаментальна в $\W$.

Пусть $q=\lim\limits_{n\to\infty}{q_n}$ в $\W$, тогда при \(z\in W_2^1[0,1]\)
$$
\langle q,z\rangle=\lim\limits_{n\to\infty}
		\int\limits_0^1 q_nz\,dx = n\int\limits_{\zeta-1/n}^{\zeta}z dx = z(\zeta).
$$
Отсюда следует, что $q=\delta_{\zeta}$.

\bigskip
\textbf{Утверждение 2.}
{\it Пусть $\gamma<1$, $\rho\in\W$, $\rho$ --- равномерно положительная функция. Тогда найдутся последовательность $\{q_n\}$ функций из \(\Gamma_\gamma\) и последовательность $\{\kappa_n\}$, где $\kappa_n\in(0,1)$ со свойством $\rho=\lim\limits_{n\to\infty}{\kappa_nq_n}$ в $\W$.}

\bigskip
\textbf{Доказательство утверждения 2.}
\medskip

Пусть $\nu\in(\gamma^+,1)$. Найдется такая постоянная функция $r\in(0,1)$, для которой выполняется неравенство $\rho(x)-r>0$ при всех $x$. Согласно утверждению 1, для любого $\eps>0$ найдется последовательность $\{\psi_n\}$ неотрицательных функций из $L_\infty[0,1]$ со свойствами $\|\psi_n\|_\nu<\eps$ при достаточно больших $n$, $\rho-r=\lim\limits_{n\to\infty}{\psi_n}$ в $\W$.
Пусть $f_n = \psi_n+r$, тогда $f^{-1}\in L_\infty[0,1]$, и $\|f_n\|_\nu\in(0,1)$, если $r$ и $\eps$ выбрать достаточно малыми.
Из утверждения 0 следует, что $\|f_n\|_\gamma\in(0,1)$. Пусть $q_n = f_n/\|f_n\|_\gamma$, тогда $q_n\in A_\gamma$ и $\|f_n\|_\gamma\cdot q_n=f_n \to \rho$.

\bigskip
\textbf{Теорема 1.}
{\it Если $\gamma<1$, то $m_\gamma=-\infty$.}

\bigskip
\textbf{Доказательство теоремы 1.}
\medskip

Будем рассматривать $\rho^*={\rm const}$ как функцию из $\W$.
Согласно утверждению 2, найдутся такие последовательности функций $q^*_n\in\Gamma_\gamma$ и чисел $\kappa_n\in(0,1)$, для которых верно $\kappa_nq^*_n\to\rho^*$.
Функция $\lambda_1(q)$ убывает при возрастании $q$, следовательно,
$\lim\limits_{n\to\infty}{\lambda_1(q^*_n)} \leq \lim\limits_{n\to\infty}{\lambda_1(\kappa_nq^*_n)} = \lambda_1(\rho^*)$, т.е.
для любого $\varepsilon>0$ при достаточно большом $n$ выполняется неравенство $\lambda_1(q^*_n) < \lambda_1(\rho^*)+\varepsilon$.


Т.к. $\rho^*$ можно выбрать сколь угодно большим, то из определения $m_\gamma$  следует, что
 $$
m_\gamma = -\infty.
 $$



\bigskip
\textbf{Утверждение 3.}
{\it Если $\gamma>1$, то \(0\in\Gamma_\gamma\).}

\bigskip
\textbf{Доказательство утверждения 3.}
\medskip

Пусть $q = \lim\limits_{n\to\infty}{q_n}$ в пространстве $\W$, где
$$
q_n(x) = \left\{
\begin{array}{ll}
n^{1/\gamma}, & x\in(0,1/n), \\
0, & \mbox{ иначе}.
\end{array}
\right.
$$
Найдем $q$:
$$
\langle q,z\rangle=\lim\limits_{n\to\infty}
		\int_0^1 q_nz\,dx = 0,
$$
т.е. $q\equiv 0$.

Т.к. $q_n\in A_\gamma$, то $q\in \Gamma_\gamma$.

\bigskip
\textbf{Теорема 2.}
{\it Если $\gamma>1$, то $M_\gamma$ --- это первое собственное значение задачи $y''+\lambda y=0$, \eqref{Karulina_eq2}.}

\bigskip
\textbf{Доказательство теоремы 1.}
\medskip

Из утверждения 3 и убывания функции $\lambda_1$ следует, что
$$
M_\gamma =\sup\limits_{q\in \Gamma_{\gamma}}\lambda_{1}(q) = \lambda_1(0).
$$


\end{document}